\documentclass[12pt]{amsart}

\usepackage{latexsym}
\usepackage{amsmath}
\usepackage{amssymb}
\usepackage{amscd}
\usepackage{multicol}
\usepackage{mathrsfs}
\usepackage{xy}
\xyoption{all}

\newtheorem{theorem}{\bf Theorem}[section]
\newtheorem{corollary}[theorem]{\bf Corollary}
\newtheorem{proposition}[theorem]{\bf Proposition}
\newtheorem{lemma}[theorem]{\bf Lemma}
\newtheorem{definition}[theorem]{\rm Definition}
\newtheorem{definition-theorem}[theorem]{\bf Theorem-Definition}

\def\bR{\mathbb{R}}
\def\bC{\mathbb{C}}
\def\bZ{\mathbb{Z}}

\def\quott({/\! /}

\def\t{\frak{t}}

\def\a{{\rm alg}}

\def\C{C_{\lambda}}

\def\z{\frak{\zeta}}

\include{include}

\begin{document}
\address{
Department of Mathematics and Statistics, University of
Regina,  Regina SK, S4S 0A2 Canada}
\email{mareal@math.uregina.ca}

\noindent {\Large \bf Connectedness of  levels
for moment maps on various classes of loop groups\footnote{
{\it MS Classification:} 53D20, 22E67}}
\bigskip\\
{\bf Augustin-Liviu Mare } \\
\smallskip \\
{\it Dedicated to Jost-Hinrich Eschenburg on his sixtieth birthday}
\title{}

\vspace{-0.5cm}

\begin{abstract}
The space $\Omega(G)$ of all based  loops in a compact  simply connected Lie group $G$ has an action of the maximal torus $T\subset G$ (by pointwise conjugation) and of the circle $S^1$
(by rotation of loops). Let $\mu : \Omega(G)\to (\t\times i\bR)^*$ be a moment map of the resulting  $T\times S^1$ action. We show that all levels (that is, pre-images of points) of $\mu$ are connected
subspaces of $\Omega(G)$ (or empty). The result holds if in the definition of 
$\Omega(G)$ loops are of class $C^{\infty}$ or of any Sobolev class $H^s$, with $s\ge 1$ (for loops of class $H^1$ connectedness of regular levels has been proved
  by Harada, Holm, Jeffrey, and the author in ~\cite{Ha-Ho-Je-Ma}).
\end{abstract}

\maketitle

\section{Introduction}
\label{intro}

Let $G$ be a compact  simply connected Lie group and
$T\subset G$ a maximal torus.
The based loop group of $G$ is the  space $\Omega(G)$ consisting of all smooth maps $\gamma: S^1\to G$ with $\gamma(1)=e$. The assignments
$$T\times \Omega(G) \to \Omega(G), \ (t,\gamma)\mapsto [S^1\ni 
z\mapsto t\gamma(z)t^{-1}]$$
and
$$S^1\times \Omega(G) \to \Omega(G), \ 
(e^{i\theta},\gamma) \mapsto [S^1\ni z\mapsto
\gamma(ze^{i\theta})\gamma(e^{i\theta})^{-1}]$$
define an action of $ T\times S^1$ on $\Omega(G)$. In fact, the latter
space is an infinite dimensional smooth symplectic manifold and the action of 
$ T\times S^1$ is Hamiltonian. Let
$$\mu: \Omega(G) \to (\t \oplus i\bR)^*$$
denote a moment map, where $\t:={\rm Lie}(T)$ and
$i\bR={\rm Lie}(S^1)$.
 Atiyah and Pressley ~\cite{At-Pr} extended the celebrated convexity theorem of Atiyah and Guillemin-Sternberg and showed that the image of
$\mu$ is the convex hull of its singular values.
Their proof's idea is to determine first the image under $\mu$ of the subspace
$\Omega_{\a}(G)\subset \Omega(G)$ whose elements are restrictions of algebraic maps from $\bC^*$ to the complexification $G^{\bC}$ of $G$:
they notice that $\mu(\Omega_{\a}(G))$ is a closed subspace of $(\t\oplus i\bR)^*$; since $\Omega_{\a}(G)$ is dense in $\Omega(G)$ (by a theorem of Segal),
they deduce from the continuity of $\mu$  that
$$\mu(\Omega(G))=\mu(\Omega_{\a}(G)).$$

The goal of this paper is to extend to $\Omega(G)$ the well-known result
which says that all levels of the moment map arising from a Hamiltonian torus action on a compact symplectic manifold are connected. That is, we will prove the following theorem.

\begin{theorem}\label{maintheorem} For any $a\in \mu(\Omega(G))$, the pre-image
$\mu^{-1}(a)$   is a connected topological subspace of $\Omega(G)$.
\end{theorem}



\noindent {\rm Remarks.}
1.  A version of Theorem \ref{maintheorem} has been proved in ~\cite{Ha-Ho-Je-Ma}. More specifically, instead of $\Omega(G)$ the authors consider there the space $\Omega_1(G)$ of all loops $S^1\to G$ of Sobolev class $H^1$. They prove that all {\it regular} levels of $\mu: \Omega_1(G) \to (\t\oplus i\bR)^*$ are  connected. 
It is not  obvious how to adapt that proof for singular levels or/and for loops of class $C^{\infty}$.  

2. One can easily see that our proof of Theorem \ref{maintheorem} works 
for $\mu: \Omega_1(G) \to (\t\oplus i\bR)^*$ (and even loops of Sobolev class $H^s$,
with
 $s\ge 1$) as well. In other words, we can prove that {\it all} levels of
$\mu : \Omega_1(G)\to (\t\oplus i\bR)^*$ are connected topological subspaces of
$\Omega_1(G)$. We decided to deal here with $\Omega(G)$ (smooth loops) rather than  
$\Omega_1(G)$ because the former  is discussed in detail in our main reference   
\cite{Pr-Se}, and the reader can make the connections directly.

We will give here an outline of the paper.
 In section \ref{two} we present basic notions and results concerning loop groups.
In section \ref{three} we define the key ingredient of the proof of Theorem \ref{maintheorem}.
This is a certain Geometric Invariant Theory (shortly G.I.T.)
quotient of $\Omega(G)$ with respect to the complexification $T^{\bC}\times \bC^*$
of $T\times S^1$. 
To define this quotient, we face the difficulty that the $ S^1$ action on
$\Omega(G)$ mentioned above does not extend to a $\bC^*$ action
(only the $T$ action extends canonically to a $T^{\bC}$ action).
However, for any $\gamma \in \Omega(G)$ there is a natural way to define the loop
$u\gamma$ for any $u \in \bC$ which is contained in the exterior of a disk with center at 0 and radius smaller than 1 (which depends on $\gamma$);
if $|u|=1$ then $u\gamma$ is given by the $S^1$ action on $\Omega(G)$ defined above.
By putting $\gamma \sim gu\gamma$, where $u$ is as before and $g\in T^{\bC}$ 
arbitrary, we obtain  an equivalence relation $\sim$ on $\Omega(G)$. The G.I.T. quotient  mentioned before is
$A/\sim$, where $A$ consists of all elements of $\Omega(G)$ which are equivalent to elements of $\mu^{-1}(a)$.  The main result of section \ref{three} is Proposition \ref{lemmaone}, which says that the natural map $\mu^{-1}(a)/(T\times S^1) \to A/\sim$ is bijective
(the  idea of the proof  belongs to Kirwan, see \cite[chapter 7]{Ki}).
In section 4 we note that the image of $(\mu^{-1}(a)\cap \Omega_{\a}(G))/(T\times S^1)$
under the map above is $(A\cap \Omega_{\a}(G))/\sim$. The former space is
connected (by a result of \cite{Ha-Ho-Je-Ma}) and we prove that the latter  is dense in 
$A/\sim$
(see Proposition \ref{propsecond}). Consequently,  $A/\sim$ is a connected topological subspace of
$\Omega(G)/\sim$. We deduce that $\mu^{-1}(a)/T\times S^1$ is connected.
Hence $\mu^{-1}(a)$ is connected as well.

 \noindent {Acknowledgements.}   I would like to  thank Jost Eschenburg for discussions about the topics of this paper. I am also grateful to the referee for carefully reading  the
 manuscript and suggesting many improvements.

 \section{Notions of loop groups}\label{two}

In this section we collect results about loop groups which will be needed later.
The details can be found in  Pressley and Segal ~\cite{Pr-Se} and/or Atiyah and Pressley
~\cite{At-Pr}.

Like in the introduction, $G$ is a compact semisimple simply connected Lie group.
We denote by $L(G)$ the space of all smooth maps $S^1\to G$ (call them loops). 
The obvious multiplication makes it into a Lie group. By $\Omega(G)$ we denote the space of all loops which map $1\in S^1$  to the unit $e$ of $G$.  
It can be naturally identified with the homogeneous space $L(G)/G$.
In fact, the presentation of $\Omega(G)$ which is most appropriate for our goals
 is
\begin{equation}\label{iden}\Omega(G)=L(G^{\bC})/L^+(G^{\bC}).\end{equation}
Here $G^{\bC}$ is the complexification of $G$ and $L(G^{\bC})$ the set of all (smooth)  loops $\alpha:S^1 \to G^{\bC}$; by  
$L^+(G^{\bC})$ we denote the   subgroup of $L(G^{\bC})$ consisting of all
$\alpha$ as above which   extend 
holomorphically for $ |\z|\le1$  (this notion is explained in detail at the beginning of the next section). Since $L(G^{\bC})$ is a complex Lie group and $L^+(G^{\bC})$ a complex Lie subgroup, equation (\ref{iden}) shows that the manifold $\Omega(G)$ has a complex structure. More precisely,   the complex structure $J_{x}$ at a point $x \in \Omega(G)$
is induced by the multiplication by $i$ in the  tangent space $T_{\alpha}L(G^{\bC})$,
where $\alpha \in L(G^{\bC})$ is such that $x= \alpha L^+(G^{\bC})$.

Let us embed $G$ into some special unitary group $SU(N)$.
We consider the 
Hilbert space $H:=L^2(S^1,\bC^N)$ and   the corresponding
 ``Grassmannian" $Gr(H)$. The latter consists of all closed vector subspaces
 of $H$ which satisfy certain supplementary properties; it turns out that
 $Gr(H)$ can be equipped with a K\"ahler (Hilbert) manifold structure (the details can be found
 in {\cite[chapter \,7]{Pr-Se}}). 
An important subspace of $Gr(H)$ is $Gr_0(H)$.
For the goals of our paper it is sufficient to mention that
$Gr_0(H)$ contains $H_+$, which is the closed vector subspace of $H$ spanned by $S^1\ni z\mapsto z^kv$, with $k\ge 0$ and $v\in \bC^N$. Also,
the   connected component of  $H_+$ in $Gr_0(H)$ consists of all vector subspaces 
$W$ of $H$ for which there exists $n\ge 0$ such that
$$z^nH_+\subset W \subset z^{-n}H_+$$
and 
$$\dim[(z^{-n}H_+)/W] =
\dim[W/(z^nH_+)].$$
In other words, if ${\mathcal G}_n$ denotes the subspace of all $W$ which satisfy
the last two equations, then the connected component of $H_+$ in  $Gr_0(H)$ is
$\bigcup_{n\ge 0} {\mathcal G}_n.$
It is important to note that via the map
$${\mathcal G}_n\ni W\mapsto W/{z^n}H_+,$$ 
the space ${\mathcal G}_n$ can be identified with the Grassmannian
$Gr_{nN}(\bC^{2nN})$ of all vector subspaces of dimension $nN$ in
\begin{equation}\label{bce}\bC^{2nN}=z^{-n}H_+/z^nH_+.\end{equation}
Also note that we have the chain of inclusions
\begin{equation}\label{filt} {\mathcal G}_0\subset  {\mathcal G}_1 \subset {\mathcal G}_2\subset \ldots.\end{equation}
Less obvious is the fact that for any $n\ge 0$, the canonical symplectic structure on 
the Grassmannian ${\mathcal G}_n$ makes this space into a symplectic submanifold of $Gr(H)$. 
The role of the above construction is revealed by the following result.
 
 \begin{proposition}\label{theom} (a) The map \begin{equation}\label{emb}\Omega(G) \to Gr(H), \  \gamma \mapsto \gamma H_+\end{equation} is an embedding, which induces on 
 $\Omega(G)$  a structure of symplectic manifold. Together with the complex structure $J$ defined above, this makes $\Omega(G)$ into a K\"ahler manifold.
 
 (b) The image of  $\Omega_{\a}(G)$ (see the introduction) under the embedding
 (\ref{emb}) is contained in $\bigcup_{n\ge 0}{\mathcal G}_n$.   
 \end{proposition} 

Based on point (b), we identify $\Omega_{\a}(G)$ with a subspace of 
$\bigcup_{n\ge 0}{\mathcal G}_n$. The inclusions  (\ref{filt}) induce the filtration
$$\Omega_{\a}(G)=\bigcup_{n\ge 0} \Omega_n, \ \Omega_0 \subset
\Omega_1 \subset \Omega_2 \subset \ldots,$$
where
$$\Omega_n:= \Omega_{\a}(G)\cap {\mathcal G}_n.$$
The space $\Omega_n$ is a closed subvariety of the Grassmannian ${\mathcal G}_n$.
We refer to the topology on $\Omega_{\a}(G)$ induced by the filtration  above 
as the {\it direct limit topology}. There is another natural topology on $\Omega_{\a}(G)$,
namely the subspace topology, induced by the inclusion
$\Omega_{\a}(G) \subset \Omega(G)$. 

The following proposition can be proved with the same arguments as 
Proposition 2.1 of  \cite{Ha-Ho-Je-Ma} (the result is  also mentioned in 
\cite[section 2]{Gu-Pr}). 

\begin{proposition}\label{fine} The direct limit topology on $\Omega_{\a}(G)$ is finer
than the subspace topology.
\end{proposition}

Let us consider again the $T\times S^1$ action on $\Omega(G)$ described
at the beginning of the paper, and the corresponding moment map 
$\mu : \Omega(G) \to (\t \oplus i\bR)^*$. 
This is uniquely determined up to an additive constant,
which will be made more precise momentarily (a standard  moment map is
described explicitly  in  \cite[section  3]{At-Pr}, but we will not need that expression
here).  For the moment, we would like to deduce from Proposition \ref{fine} a
result which will be useful later. Namely, let us take $a\in \mu(\Omega(G))$;
 by \cite[Proposition 3.4]{Ha-Ho-Je-Ma},  $\mu^{-1}(a)\cap \Omega_{\a}(G)$
 is  a connected subspace of $\Omega_{\a}(G)$ with respect to the direct limit topology.
 We deduce: 
 
\begin{proposition}\label{conncor} For any $a\in \mu(\Omega_{\a}(G))$, the space 
$\mu^{-1}(a)\cap \Omega_{\a}(G)$   is a connected topological subspace of $\Omega(G)$.
\end{proposition}

There is also an action of $T\times S^1$ on each ${\mathcal G}_n$, $n\ge 0$, which can be described as follows. We fix a basis, say $b_1,\ldots,b_N$,  of $\bC^N$,
and consider the induced basis $z^{k}b_j$, $-n\le k \le n-1$, $1\le j\le N$, of
$\bC^{2nN}$ (see equation (\ref{bce})). 
 The action of  $T$ on ${\mathcal G}_n$ is induced by
 \begin{equation}\label{te}t. (z^{k}b_j):= z^{k}(tb_j),\end{equation}
 for any $t\in T$ and $k,j$ as above; the action of $S^1$ is induced by
 \begin{equation}\label{e}e^{i\theta}. (z^{k}b_j) := (e^{i\theta}z)^kb_j=z^k e^{ik\theta}b_j\end{equation}
 for all $e^{i\theta}\in S^1$. 
This $T\times S^1$ action is the restriction of 
an obvious $T^{\bC}\times \bC^*$  action: namely, in equation (\ref{te}) we take
$t\in T^{\bC}$ and in equation (\ref{e}) we replace $e^{i\theta}$ by an arbitrary element of $\bC^*$. The $T^{\bC}\times \bC^*$ action turns out to be linear with respect to the 
Pl\"ucker embedding
of ${\mathcal G}_n$ (see {\cite[section  \,4]{At-Pr}}). Thus, the $T\times S^1$ action is Hamiltonian.
We pick
$$\mu_n : {\mathcal G}_n \to (\t \oplus i\bR)^*$$ a moment map, which is
again uniquely determined up to an additive constant.
We can arrange the constants in such a way that
if $m<n$ then
$$\mu_n|_{{\mathcal G}_m}=\mu_m.$$
The reason is that ${\mathcal G}_m$ is a $T\times S^1$-invariant symplectic
submanifold of ${\mathcal G}_n$.
We  obtain the map
$\tilde{\mu}: \bigcup_{n\ge 0}{\mathcal G}_n\to (\t \oplus i\bR)^*$
such that $\tilde{\mu}|_{{\mathcal G}_n}=\mu_n$, for all $n\ge 0$. 
The map $\tilde{\mu}$  is uniquely determined up to an additive constant.
The following proposition relates the moment maps $\mu$ and $\tilde{\mu}$.

\begin{proposition}\label{muuu} We can choose $\mu$ and $\tilde{\mu}$ such that
$$\mu|_{\Omega_{\a}(G)}= \tilde{\mu}|_{\Omega_{\a}(G)}.$$
\end{proposition} 

\begin{proof} The idea of the proof is that there exists a submanifold
$Gr_{\infty}(H)$ of $Gr(H)$ acted on smoothly by $T\times S^1$ and such that
\begin{itemize}
\item $Gr_0(H)\subset Gr_{\infty}(H)$ and the inclusion is $T\times S^1$ equivariant
\item there exists  $\hat{\mu}: Gr_{\infty}(H)\to \t\oplus i\bR$ which is a moment map for the $T\times S^1$ action
\item the image of $\Omega(G)$ under the inclusion (\ref{emb}) is contained in
$Gr_{\infty}(H)$.
\end{itemize}
It is worth noticing that $Gr(H)$ does not admit a smooth action of $T\times S^1$;
only its subspace $Gr_{\infty}(H)$ does (see \cite[section 7.6]{Pr-Se}). This is why we need to use the latter space
in our proof.
 
We deduce that $\hat{\mu}|_{\Omega(G)}$ differs from $\mu$ by a constant;
the same can be said about $\hat{\mu}|_{{\mathcal G}_n}$ and $\mu_n$, for any 
$n\ge 0$.
The result follows.
\end{proof} 

Let us consider again the $T^{\bC}\times \bC^*$ action on ${\mathcal G}_n$
defined above. Any of the inclusions ${\mathcal G}_n\subset {\mathcal G}_{n+1}$
is equivariant. Thus, we have an action of $T^{\bC}\times \bC^*$ on $\bigcup_{n\ge 0}{\mathcal G}_n$. 
The same group acts  on $\Omega_{\a}(G)$, as follows.
We take into account that
\begin{equation}\label{om}\Omega_{\a}(G)=L_{\a}(G^{\bC})/L_{\a}^+(G^{\bC})
\end{equation}
where $L_{\a}(G^{\bC})$ is the space of all  algebraic maps
$\alpha:\bC^*\to G^{\bC}$ and $L_{\a}^+(G^{\bC})$ the subgroup consisting
of those $\alpha$ which can be extended holomorphically to $\bC$. Then the action we are
referring to is
\begin{equation}\label{acti}T^{\bC}\times \bC^*\times \Omega_{\a}(G)\ni(g,u,\alpha L_{\a}^+(G^{\bC}))
\mapsto [\bC^*\ni \z \mapsto g\alpha(u\z)]L_{\a}^+(G^{\bC}).\end{equation}
The following result will be needed later.

\begin{proposition} The inclusion $\Omega_{\a}(G)\subset
\bigcup_{n\ge 0}{\mathcal G}_n$ defined in Proposition \ref{theom} (b) is
$T^{\bC}\times \bC^*$ equivariant.
\end{proposition}

\begin{proof} Take $\gamma \in \Omega_{\a}(G)$, which is  of the form
$$S^1\ni z\mapsto \gamma(z)=\sum_{-k_0\le k\le k_0} A_k z^k,$$
where $k_0\ge 0$.
Here  $A_k$ are $N\times N$ matrices with entries in $\bC$. 
The subspace $\gamma H_+$ of $H$ has the property that 
$$z^{n}H_+ \subset \gamma H_+\subset z^{-n} H_+,$$
for some $n \ge 0$. Any element $v$ of $H_+$ has a Fourier expansion of the
form $v=\sum_{m\ge 0}v_m z^m$, where $v_m \in \bC^N$, for all $m\ge 0$.
Then $$\gamma v = \sum_{m\ge 0, k\in \bZ}(A_kv_m)z^{k+m}.$$
The corresponding element of $(\gamma H_+)/z^n H_+$ is
 $$[\gamma v]=\gamma v {\rm mod} z^n H_+=
 \sum_{m \ge 0, -k_0\le k\le k_0, k+m\le n-1}
 (A_kv_m)z^{k+m}.$$
 In this sum we have $m\le n-k-1\le n+k_0-1$.
 Thus, to describe all elements of $\gamma H_+/z^m H_+$, it is sufficient to take
 $v$ of the form $$v=\sum_{0\le m \le n+k_0-1} v_m z^m.$$ 
  If $t\in T^{\bC}$, then
 $$t. [\gamma v] =\sum_{m\ge 0, -k_0\le k\le k_0, k+m \le n-1}t(A_kv_m)z^{k+m} =\sum_{m\ge 0, -k_0\le k\le k_0, k+m \le n-1}(tA_k)v_mz^{k+m}
=[(t.\gamma) v].$$
Consequently, $t.(\gamma H_+) =(t.\gamma) H_+$. 
If $u\in \bC^*$, then
$$u.[\gamma v] =\sum_{m\ge 0, -k_0\le k\le k_0, k+m \le n-1}(A_kv_m)(uz)^{k+m}=\sum_{m\ge 0, -k_0\le k\le k_0, k+m \le n-1}(u^kA_ku^mv_m)z^{k+m}.
$$
This is the same as $[(u.\gamma) \tilde{v}]$, where 
$$\tilde{v}=\sum_{0\le m \le n+k_0-1} u^mv_m z^m.$$
Consequently, $u.(\gamma H_+) = (u.\gamma)H_+$.
\end{proof}

Finally, let us pick $B\subset G^{\bC}$ a Borel subgroup with $T\subset B$.
The presentation (\ref{om}) of $\Omega_{\a}(G)$ allows us to define on the latter space a natural action of the group 
$${\mathcal B}_+:=\{\alpha\in L^+_{\a}(G^{\bC}) \ : \ \alpha(0) \in B\}$$ on $\Omega_{\a}(G)$.
The orbit decomposition is
$$\Omega_{\a}(G) = \bigcup_{\lambda\in \check{T}}C_{\lambda},$$
where the union is disjoint and $$C_{\lambda}:= {\mathcal B_+} \lambda$$ is called a 
{\it Bruhat cell}. Here $\check{T}$ denotes the lattice of group homomorphisms
$S^1\to T$. 
The space $\C$  is really a (finite dimensional) cell, being homeomorphic to $\bC^{r}$ for some $r$. In this paper, by $\overline{\C}$ we will always mean  the
closure of $\C$ in the direct limit topology (see above). The following property of the Bruhat cells will be needed later. 

\begin{proposition}\label{subva} For any $\lambda\in \check{T}$, there exists $n\ge 1$ such that
$\overline{\C}$ is contained in ${\mathcal G}_n$ as a $T^{\bC}\times \bC^*$-invariant closed subvariety. 
\end{proposition}

This can be proved as follows. There exists $n\ge 1$ such that
$\C \subset {\mathcal G}_n$, because ${\mathcal B}_+$ leaves each $\Omega_k$,
$k\ge 0$, 
invariant (see \cite[Lemma   \,3.3.2]{Ko}).
The space $\C$ is a locally Zariski closed subspace of $\Omega_{\a}(G)$
(see \cite[Proposition 2.13 and Theorem 3.1]{Mi}), thus also of $\Omega_n$ and
of ${\mathcal G}_n$. Consequently, the  closures of $\C$ in the Zariski, respectively differential topology of ${\mathcal G}_n$ are equal.

Another  result concerning the Bruhat cells  is the following proposition
(cf. ~\cite[section 1]{At-Pr}, see also {\cite[proof of Proposition 3.4]{Ha-Ho-Je-Ma}}).

\begin{proposition}\label{last} For any $\lambda_1,\lambda_2\in \check{T}$ there exists
$\lambda\in \check{T}$ such that
$$C_{\lambda_1}\subset \overline{\C} \ {\rm and } \  C_{\lambda_2}\subset 
\overline{\C}.$$
Consequently, for any $x,y\in\Omega_{\a}(G)$ there exists $\lambda\in \check{T}$
such that both $x$ and $y$ are in $\overline{\C}$.
\end{proposition}










 \section{The equivalence relation $\sim$}\label{three}

We begin with the following definition. Take $0<r\le 1$. We say that  a 
free loop $S^1\to G^{\bC}$ {\it extends holomorphically for $|\z|\ge r$ } if 
it is the restriction of a map 
$${\alpha} : \{\z\in \bC\cup\infty  \ : \ |\z| \ge r\} \to G^{\bC}$$
which is continuous,  holomorphic on $ \{\z\in \bC \cup\infty  \ : \ |\z| > r\}$ 
and smooth on $\{\z \in \bC \ : \ |\z| =r\}$;
the same terminology is adopted  if we take $r\ge 1$ and replace ``$\ge$" and ``$ >$" by ``$\le$", respectively ``$<$" (and also $\bC\cup\infty$ by $\bC$).

Let $L^-(G^{\bC})$ denote the subspace of $L(G^{\bC})$ consisting of those
$\alpha$  which extend holomorphically for $|\z|\ge 1$ in the sense of the definition above. 
One knows that any  $\alpha \in L(G^{\bC})$ can be written as
$$\alpha = \alpha_- \lambda \alpha_+,$$
where $\alpha_- \in L^-(G^{\bC})$,  $\alpha_+ \in L^+(G^{\bC})$, and
$\lambda$ is a group homomorphism $S^1 \to T$ (see
{\cite[Theorem 8.1.2]{Pr-Se}}).
By using the presentation (\ref{iden}), the elements of $\Omega(G)$ are cosets of the form
$\alpha_-\lambda  L^+(G^{\bC})$, where $\alpha_-$ and $\lambda$ are as above.
The following lemma will be used later.

\begin{lemma}\label{lemmaoone} Take $\alpha_-,\beta_-$  in $L^-(G^{\bC})$ and 
$\lambda,\mu :S^1\to T$  group homomorphisms such that
$$\alpha_-\lambda L^+(G^{\bC}) = \beta_-\mu L^+(G^{\bC}).$$
Let $r$ be a  strictly positive real number.

(a) Assume that $r<1$. 
If $\alpha_-$ extends holomorphically for $|\z|\ge r$ then 
$\beta_-$ extends holomorphically for $|\z|\ge r$ as well.

(b)  Assume that $r\ge 1$ or $r<1$ and $\alpha_-$ extends holomorphically for 
$|\z|\ge r$. 
For any $u\in \bC^*$ with $|u|\ge r$ we have
$$[S^1\ni z\mapsto \alpha_-(uz)]\lambda L^+(G^{\bC}) =
[S^1\ni z\mapsto  \beta_-(uz)]\mu L^+(G^{\bC}).$$
\end{lemma}

\begin{proof} 
We have
\begin{equation}\label{alp}\alpha_-\lambda  = \beta_-\mu \alpha_+,
\end{equation}
where $\alpha_+\in L^+(G^{\bC})$. 

(a) The loops $\lambda$ and $\mu$ are one-parameter subgroups in $T$, thus they have 
obvious (holomorphic) extensions  to group homomorphisms $\bC^*\to T^{\bC}$.
From (\ref{alp}) we deduce that
$\beta_-$ is the restriction of a function holomorphic on the annulus
$$\{\z \in \bC  \ : \ r<|\z| <1\}$$ and continuous on
the closure of this space. Consequently, the map
$\xi\mapsto {\beta}_-(\frac{1}{\xi})$ extends holomorphically for
$|\xi|\le \frac{1}{r}$,
that is, $\beta_-$  extends holomorphically   
for $|\z| \ge r$. Indeed, let us consider again the embedding 
$G\subset SU(N)$, as in section \ref{two}.  The resulting embedding 
$G^{\bC}\subset {\rm Mat}^{N\times N}(\bC)$ is holomorphic. We  use the
following claim:

\noindent{\it Claim.} If $f:\{\xi\in \bC \ : \ |\xi|\le \frac{1}{r}\}\to \bC$ is a continuous
function which is holomorphic on $\{\xi\in \bC \ : \ |\xi|<\frac{1}{r}, |\xi|\neq 1\}$,
then $f$ is holomorphic on $\{\xi\in \bC \ : \ |\xi|<\frac{1}{r}\}$.

This can be proved by comparing the Laurent series of $f$ on
$\{\xi \in \bC \ : \ |\xi|<1\}$, respectively $\{\xi \in \bC \ : \ 1<|\xi|<\frac{1}{r}\}$.
The series are equal, since the  coefficients of both of them are equal to
$\frac{1}{2\pi i}\int_{|\xi|=1}\frac{f(\xi)}{\xi^k} d\xi$, $k\in \bZ$ (by a uniform continuity argument). Thus, the radius of convergence of the first of the two series (which is actually the Taylor series of $f$ 
around $0$)
is at least equal to $\frac{1}{r}$. The claim is proved. 

(b) From equation (\ref{alp}) we deduce that
$\alpha_+$ extends holomorphically to $\bC$.
The reason is that the entries of the $N\times N$ matrix
$\alpha_+=\mu^{-1}\beta_-^{-1}\alpha_-\lambda$ are $\bC$-valued functions 
which are continuous on $\bC$ and holomorphic on
$\bC\setminus\{\xi \in \bC \ : \ |\xi|=1\}$; by the same argument as in the claim above, they are
holomorphic on the whole $\bC$.
Again from equation (\ref{alp}), we deduce that
$$\alpha_-(uz)\lambda(uz)=\beta_-(uz)\mu(uz)\alpha_+(uz),$$
for all $z\in S^1$. The map $S^1\ni z\mapsto \alpha_+(uz)$ is in $L^+(G^{\bC})$.
We only need to notice that $$\lambda(uz)=\lambda(z)\lambda(u), \ \mu(uz)=\mu(z)\mu(u).$$
\end{proof}

\begin{definition}\label{defi}  (a) Take $x\in \Omega(G)$ and $u\in \bC^*$.
We say that the pair $(u,x)$ is {\rm admissible} if 
\begin{itemize} 
\item $|u|\ge 1$
\end{itemize}
or
\begin{itemize}
\item 
$|u|< 1$ and  $x=\alpha_-\lambda L^+(G^{\bC})$, where
$\alpha_-\in L^-(G^{\bC})$ extends holomorphically for $|\z|\ge |u|$ and 
$\lambda:S^1\to T$ is a group homomorphism.
\end{itemize}
 If $(u,x)$ is as above and $g\in T^{\bC}$, then
$$gux:=g[S^1\ni z\mapsto \alpha_-(uz)]\lambda L^+(G^{\bC})$$
is an element of $\Omega(G).$

(b) Take $x,y\in \Omega(G)$. We say that 
$$x\sim y$$ if there exist  $u\in \bC^*$ and $g\in T^{\bC}$ such that
$(u,x)$ is an admissible pair and
$y=gux$.
\end{definition}

\noindent{\bf Remark.} We can also express $gux$ as
$$gux:=g[S^1\ni z\mapsto (\alpha_-\lambda)(uz)] L^+(G^{\bC}),$$
because $\lambda$ is a group homomorphism $\bC^*\to T^{\bC}$.

Note that by Lemma \ref{lemmaoone}, the definition of $gux$ in part (a)  is independent of the choice of the representative $\alpha_-\lambda$  of $x\in L(G^{\bC})/L^+(G^{\bC})$.  The following lemma shows that $\sim$ is an equivalence relation.

\begin{lemma}\label{ifx} (a) If $x\in \Omega(G)$,  $u\in \bC^*$ and
$g\in T^{\bC}$
such that $(u,x)$ is admissible, then $(u^{-1},gux)$ is admissible and we have
$$g^{-1}u^{-1}(gux)=x.$$

(b) If $x\in \Omega(G)$, $u_1,u_2\in \bC^*$, and $g_1,g_2\in T^{\bC}$ such that
$(u_1,x)$  and $(u_2, g_1u_1x)$ are admissible, then
$(u_1u_2,x)$ is admissible and
$$(g_1g_2)(u_1u_2)x=g_2u_2(g_1u_1x).$$
\end{lemma} 

\begin{proof} 
(a) We can assume that $g=1$.  We write $x=\alpha_-\lambda L^+(G^{\bC})$.
Assume first that $|u|\ge 1$.  
The loop
$S^1\ni z\mapsto \alpha_-(uz)$  extends holomorphically for
$|\z|\ge \frac{1}{|u|}$ by $\z \mapsto {\alpha}_-(u\z)$.
The case  $|u|< 1$ is even easier to analyze. 
Verifying that $u^{-1}(ux)=x$ is equally easy.

(b) We can assume that $g_1=g_2=1$. Again we write $x=\alpha_-\lambda L^+(G^{\bC})$. 
It is sufficient to analyze the case when $| u_1u_2|<1$.
Thus, at least one of the numbers $|u_1|$ and $|u_2|$ is strictly less than 1.
We distinguish the following two cases.

\noindent {\it Case 1.} $|u_2|<1$. The loop $S^1\ni z \mapsto \alpha_-(u_1z)$ is well defined 
and extends holomorphically  for $|\xi|\ge |u_2|$. 
Let $\tilde{\alpha}:\{\xi\in \bC\cup \infty \ : \ |\xi|\ge |u_2|\} \to G^{\bC}$ be an extension of this loop. 
The map $\hat{\alpha} : \{\z\in \bC \cup \infty \ : \ |\z|\ge |u_1u_2|\} \to G^{\bC}$
given by
$$\hat{\alpha}(\z)
=
\begin{cases}
\tilde{\alpha}(u_1^{-1}\z), & {\rm if} \  |u_1u_2|\le |\z|\le |u_1|\\
\alpha_-(\z), & {\rm if} \ |u_1|\le |\z|
\end{cases}
$$
is the desired extension of $\alpha_-$ for
$ |\z| \ge |u_1u_2|$ (note that
$\hat{\alpha}$ is holomorphic on $|\z|>|u_1u_2|$, since it is 
continuous and is  holomorphic on the complement of the circle
$\{\z\in \bC \ : \ |\z|=|u_1|\}$). 

\noindent {\it Case 2.} $|u_2|\ge 1$. This implies $|u_1|<1$. We notice that the pair
$( u_1, u_2x)$ is admissible: indeed, by hypothesis, the loop $S^1\ni z\mapsto \alpha_-(u_2z)$
extends holomorphically for $|u_2\z|\ge |u_1|$, hence also for $|\z|\ge |u_1|$. 
The pair $( u_2,x)$ is admissible too. From the result proved in case 1 we deduce that
$( u_1u_2, x)$ is admissible. 

The equation
$u_2(u_1x)=(u_1u_2)x$ is straightforward.
\end{proof} 

The following result relates the equivalence relation $\sim$ to the $T\times S^1$ action on $\Omega(G)$ (see section 1).  
 
\begin{lemma}\label{twotwo} Take $\gamma \in \Omega(G)$. If  $\theta \in \bR$, then
the pair $(e^{i\theta}, \gamma)$ is admissible.
If $t\in T$, then the loop $te^{i\theta}\gamma$ given by Definition   \ref{defi} (b) 
can be expressed as 
$$te^{i\theta}\gamma =t \gamma^{\theta}t^{-1}.$$
Here the right-hand side is given by
$$(t\gamma^\theta t^{-1})(z)=t\gamma(ze^{i\theta})\gamma(e^{i\theta})^{-1}t^{-1},$$
for all $z\in S^1$. 
\end{lemma} 
\begin{proof}
There exist $\alpha_-\in L^-(G^{\bC})$ and  $\lambda :S^1\to T$ a group homorphism  such that the image of
$\gamma$ under the isomorphism (\ref{iden}) is $\alpha_-\lambda L^+(G^{\bC})$.
This means that 
$$\alpha_-\lambda = \gamma \alpha_+,$$
for some $\alpha_+\in L^+(G^{\bC})$.
We deduce that for any $z\in S^1$ we have
$$[t\alpha_-(ze^{i\theta})\lambda(z)]\lambda(e^{i\theta})
=[t\gamma^{\theta}(z)t^{-1}]t\gamma(e^{i\theta})\alpha_+(ze^{i\theta}).$$
In other words, via the isomorphism (\ref{iden}), to $t\gamma^{\theta}t^{-1}$ corresponds the coset of 
$$t[S^1\ni z\mapsto \alpha_-(ze^{i\theta})]\lambda,$$
which is the same as $te^{i\theta}(\alpha_-\lambda L^+(G^{\bC}))$.
\end{proof}

We now denote by $A$ the set of all $x\in \Omega(G)$ with $x\sim y$, for some $y\in \mu^{-1}(a)$.
We are interested in the quotient space $A/\sim$ and the (natural) map 
$\mu^{-1}(a)/(T\times S^1) \to A/\sim$ which assigns to the coset of $x\in \mu^{-1}(a)$
the equivalence class of $x$. By Lemma \ref{twotwo}, this map is  well defined.

\begin{proposition}\label{lemmaone} The natural map 
$$\mu^{-1}(a)/(T\times S^1) \to A/\sim$$
is  bijective.
\end{proposition}
\begin{proof}
Only the injectivity has to be proved. We have to show that if
$x,y\in \mu^{-1}(a)$ with $x\sim y$, then $y=te^{i\theta}x$, where 
$(t,e^{i\theta})\in T\times S^1$.
By  Definition \ref{defi} we have
$$y=gux$$
for some $u\in \bC^*$ and $g\in T^{\bC}$.
We write $g=\exp(w_1)\exp(iw_2)$ and $u=e^{i\alpha}e^{-\beta}$, where 
$w_1,w_2\in\t$ and $\alpha,\beta\in \bR$ (here we see $-\beta$ as $i(i\beta)$). 
Since the pair $(u,x)$ is admissible and  $|u|=|e^{-\beta}|$, the pair
$(e^{-\beta},x)$ is admissible too. By Lemma \ref{ifx} (b) we have
$$y=\exp(w_1)e^{i\alpha} (\exp(iw_2)e^{-\beta} x).$$
Thus, it  is sufficient to assume that
$$y=\exp(iw_2)e^{-\beta} x.$$
Moreover, without loss of generality we assume that $$\beta\ge 0,$$
because if contrary we write $x=\exp(-iw_2)e^{\beta}y$ (by Lemma \ref{ifx} (a)).
Let us consider the function $h:[0,1]\to \bR$,
$$h(s)=[\mu(\exp(isw_2)e^{-s\beta} x)-a](w_2,i\beta),$$
where $0\le s \le 1$. Notice that  $h(s)$ is well defined for any $s$
with $0\le s \le 1$: indeed, the pair $(e^{-\beta}, x)$ is admissible
hence, because $e^{-s\beta}\ge e^{-\beta}$, the pair $(e^{-s\beta},x)$
is admissible too.  Since $\mu(x)=\mu(y)$, we have $h(0)=h(1)=0$.
Consequently,  there exists
$s_0$ in the interval $(0,1)$  such that $h'(s_0)=0$. 
We  denote 
\begin{equation}\label{k} x_0:=\exp(is_0w_2)e^{-s_0\beta} x.\end{equation} 

\noindent {\it Claim.} We have
$$\frac{d}{ds}\bigg|_{s_0}\exp(isw_2)e^{-s\beta}x
= J_{x_0}((w_2,i\beta) .x_0),$$
where $J_{x_0}$ is the complex structure at $x_0$ (see section 2) and 
$$(w_2,i\beta).x_0:=\frac{d}{ds}\bigg|_0[\exp(sw_2)e^{is\beta}x_0]$$
arising from the infinitesimal action of $T\times S^1$ on $\Omega(G)$.

The claim can be proved as follows. Write $x_0=\alpha_-\lambda L^+(G^{\bC})$, where
$\alpha_-\in L^-(G^{\bC})$ and $\lambda : S^1 \to T$ is a group homomorphism.
By using Lemma \ref{ifx} and the remark following Definition \ref{defi}, we have  
$$\exp (is w_2)e^{-s\beta}x=
\exp (i(s-s_0)w_2)e^{-(s-s_0)\beta}x_0
=\exp (i(s-s_0)w_2) (\alpha_-\lambda)_{-(s-s_0)} L^+(G^{\bC}).$$
Here we have denoted 
$$(\alpha_-\lambda)_{-(s-s_0)}(z):=(\alpha_-\lambda)(e^{-(s-s_0)\beta}z)$$
for all $s$ in the interval $(0,1)$ and all $z\in S^1$.
By the definition of the complex structure $J$ (see section 2), it is sufficient to prove that 
\begin{equation}\label{eql}\frac{d}{ds}\bigg|_{s_0}[\exp (i(s-s_0)w_2) (\alpha_-\lambda)_{-(s-s_0)}]
= i \frac{d}{ds}\bigg|_{0}[\exp (sw_2) (\alpha_-\lambda)_{is}]
,\end{equation}
where 
$$(\alpha_-\lambda)_{is}(z):=(\alpha_-\lambda)(e^{is\beta}z)$$
for all $s\in \bR$ and all $z\in S^1$.
By using the Leibniz rule, the left-hand side of (\ref{eql}) is
$$\frac{d}{ds}\bigg|_{0}[\exp (isw_2) (\alpha_-\lambda)_{-s}]
=i\frac{d}{ds}\bigg|_{0}[\exp (sw_2)](\alpha_-\lambda)
+i\frac{d}{ds}\bigg|_0[(\alpha_-\lambda)_{is}].
$$
Here we have used that
$$\frac{d}{ds}\bigg|_{0}[\exp (isw_2)]
=iw_2=i\frac{d}{ds}\bigg|_{0}[\exp (sw_2)]$$
and also that
$$\frac{d}{ds}\bigg|_0(\alpha_-\lambda)(e^{-s\beta}z)
=i\frac{d}{ds}\bigg|_0(\alpha_-\lambda)(e^{is\beta}z),$$
for all $z\in S^1$ (the last equation follows from the fact that
$\alpha_-\lambda$ is holomorphic on the exterior of a closed disk with center at 0
and radius strictly smaller than 1).
The claim is proved.

From the claim we deduce as follows:
\begin{align*}h'(s_0)&=(d\mu)_{x_0}(\frac{d}{ds}\bigg|_{s_0}(\exp(isw_2)e^{-s\beta}x))(w_2,i\beta)\\
{}&=\omega_{x_0}(\frac{d}{ds}\bigg|_{s_0}(\exp(isw_2)e^{-s\beta}x),(w_2, i\beta).x_0)
\\{}&=\omega_{x_0}(J_{x_0}((w_2,i\beta) .x_0),(w_2, i\beta).x_0)
\\{}&=\langle (w_2, i\beta).x_0, (w_2, i\beta).x_0 \rangle,\end{align*}
where $\omega$ denotes the symplectic form and $\langle \ , \  \rangle$   the K\"ahler metric on $\Omega(G)$ (see Proposition 2.1).
We deduce that $$(w_2, i\beta).x_0=0$$
which, according to the claim above,   implies that
$$\frac{d}{ds}\bigg|_0\exp(isw_2)e^{-s\beta}x_0=0.$$
From this we deduce that
\begin{equation}\label{kkk}\exp(isw_2)e^{-s\beta}x_0=x_0\end{equation}
for all $ s\le 0$ (note that for any such $s$, the pair
$(e^{-s\beta},x_0)$ is admissible, since 
$e^{-s\beta}\ge 1$). Indeed, by using Lemma \ref{ifx} we deduce that for any $s_1 \le 0$ we have 
\begin{align*}\frac{d}{ds}\bigg|_{s_1}\exp(isw_2)e^{-s\beta}x_0 & =\frac{d}{ds}\bigg|_{s_1}(\exp(is_1w_2)e^{-s_1\beta})\exp(i(s-s_1)w_2)e^{-(s-s_1)\beta}x_0
\\
{} &  =d(\exp(is_1w_2)e^{-s_1\beta})_{x_0}(
\frac{d}{ds}\bigg|_0(\exp(isw_2)e^{-s\beta}x_0))\\
{} & =0.\end{align*}
Here we have used the (differential of the) map 
$\exp(is_1w_2)e^{-s_1\beta}:\Omega(G)\to \Omega(G)$ given by
$$\gamma\mapsto \exp(is_1w_2)e^{-s_1\beta}\gamma,$$
which is well defined, since $e^{-s_1\beta}\ge 1$.

By Lemma \ref{ifx} (a), equation (\ref{kkk}) implies that the pair $(e^{-s\beta},x_0)$ is admissible for any $s\ge 0$; 
moreover, equation (\ref{kkk}) holds for all $s\ge 0$ as well. We make $s=-s_0$ in (\ref{kkk}) and deduce 
$x=x_0$; then we make $s=1-s_0$ and deduce $y=x_0$. We conclude
$$x=y$$ and the proof is finished.
\end{proof}

\noindent {\bf Remark.} 
Let $M$ be a compact K\"ahler manifold acted on by a complex Lie group $G$, which is the complexification of a compact Lie group $K$, in such a way that the action of 
$K$ on $M$ is Hamiltonian. Kirwan has proved that if $x,y \in M$ have the same
image under the moment map and are on the same $G$ orbit, then they are on the same $K$-orbit (see \cite[Lemma 7.2]{Ki}). We have used above the idea of her proof.
Kirwan's result cannot be used directly in our context: first, $\Omega(G)$ is
not a compact manifold; second, and most importantly, the $T\times S^1$ action on
$\Omega(G)$ does not extend in any reasonable way to a  $T^{\bC}\times \bC^*$ action. We are substituting this action by the equivalence relation $\sim$.

 \section{Connectedness of $A/\sim$ and of $\mu^{-1}(a)$}

We start with the following proposition.

\begin{proposition}\label{propfirst} 
(a) If $x\in \Omega_{\a}(G)$ then the pair $(u,x)$ is admissible (in the sense of Definition \ref{defi})
for any 
$u\in \bC^*$. The map $$T^{\bC}\times \bC^*\times 
\Omega_{\a}(G) \to \Omega_{\a}(G), \ (g,u,x)\mapsto gux$$
is the  action of   $T^{\bC}\times \bC^*$ on $\Omega_{\a}(G)$
defined in section \ref{two} (see equation (\ref{acti})).

(b) The image of
$(\mu^{-1}(a)\cap \Omega_{\a}(G))/(T\times S^1)$ under the map
in Proposition  \ref{lemmaone} is
$(A\cap\Omega_{\a}(G))/\sim$. The latter space  is
a connected topological subspace of $\Omega(G)/\sim$.
\end{proposition}

\begin{proof} Point (a) follows  from equations (\ref{om}) and  (\ref{acti}) and
the remark following Definition \ref{defi}.
To prove the first assertion  of (b), we only need to note that if $x\in \Omega_{\a}(G)$
and $y\in \Omega(G)$ such that $x\sim y$, then $y\in \Omega_{\a}(G)$. To prove the second assertion of (b),
we note that the natural map
$$(\mu^{-1}(a)\cap \Omega_{\a}(G))/(T\times S^1)\to \Omega(G)/\sim$$
is continuous. We use Proposition \ref{conncor}.
\end{proof}

The key result of this section is
\begin{proposition}\label{propsecond} The subspace  $(A\cap\Omega_{\a}(G))/\sim$ of $A/\sim$ is dense (both spaces have the topology of subspace of $\Omega(G)/\sim$).
\end{proposition}
Combined with Proposition \ref{propfirst} (b), this implies
\begin{corollary}\label{cor} The space $A/\sim$ is a connected topological subspace
of $\Omega(G)/\sim$.
\end{corollary}

In turn, this implies the main result of the paper, as follows.

\noindent {\it Proof of Theorem \ref{maintheorem}.} 
The natural map
\begin{equation}\label{induce}\mu^{-1}(a)/(T\times S^1)\to \Omega(G)/\sim
\end{equation}
is continuous, one-to-one, and its image is $A/\sim$ (by Proposition \ref{lemmaone}).
Since $A/\sim$ is connected (see the previous corollary), we deduce 
that   $\mu^{-1}(a)/(T\times S^1)$ is  connected as well.
  Consequently, $\mu^{-1}(a)$ is a connected topological subspace of
$\Omega(G)$. \hfill $\square$

The rest of the section is devoted to the proof of Proposition \ref{propsecond}. 
First, if $\lambda \in \check{T}$, we say that a point $x\in \overline{\C}$ is
$(\mu-a)$-{\it semistable} if
\begin{equation}\label{ecu1}\overline{(T^{\bC}\times \bC^*)x} \cap (\mu^{-1}(a)\cap 
\overline{\C})\neq \emptyset.\end{equation}
Here the closure is taken in $\Omega_{\a}(G)$ with respect to the direct limit topology.  
We may assume  that $\overline{\C}$ is contained in the Grassmannian 
${\mathcal G}_n$  as a 
$T^{\bC}\times
\bC^*$-invariant closed subvariety (see Proposition \ref{subva}).  Then $x$ is $(\mu-a)$-semistable if and only
if it is $(\mu_n-a)$-semistable in the usual sense, that is, if
$$\overline{(T^{\bC}\times \bC^*)x} \cap (\mu_n^{-1}(a)\cap 
\overline{\C})\neq \emptyset$$ (see for instance
{\cite[chapter 7]{Ki}}). This follows immediately from the fact that 
$\mu$ and $\mu_n$ coincide on $\overline{\C}$, by 
Proposition \ref{muuu}. We denote by $\overline{\C}^{ss}$ the set of all semistable points in
$\overline{\C}$. We also consider  the set  ${\mathcal G}_n^{ss}$ of all 
$(\mu_n-a)$-semistable points in ${\mathcal G}_n$. We have 
\begin{equation}\label{semis}\overline{\C}^{ss}=\overline{\C}\cap {\mathcal G}_n^{ss}.\end{equation}
The following description of the semistable set of $\overline{\C}$ will be needed later. 
\begin{lemma}\label{lema} We have
$$A\cap \overline{\C}=\overline{\C}^{ss}.$$
\end{lemma}

\begin{proof} 
By Proposition \ref{propfirst} (b), we have
$$A\cap \Omega_{\a}(G) =(T^{\bC}\times \bC^*)(\mu^{-1}(a)\cap
\Omega_{\a}(G)).$$
Consequently, a point $x\in\Omega(G)$ is in  $A\cap \overline{\C}$ 
if and only if $x\in [(T^{\bC}\times \bC^*)\mu^{-1}(a)]\cap \overline{\C}$.
The latter set is obviously equal to 
$(T^{\bC}\times \bC^*)(\mu^{-1}(a)\cap \overline{\C})$,
which is the same as $\overline{\C}^{ss}$ (by 
{\cite[Theorems 7.4 and 8.10]{Ki}}, applied for
the Grassmannian ${\mathcal G}_n$ which contains $\overline{\C}$ as a 
$T^{\bC}\times
\bC^*$-invariant closed subvariety, as  indicated above).
\end{proof}

We are now ready to prove Proposition \ref{propsecond}.

\noindent {\it Proof of Proposition \ref{propsecond}.} We show that in any open subset
$V$ of
$A/\sim$ there exists an element of $(A\cap\Omega_{\a}(G))/\sim$.
Since $A/\sim$ is equipped with the topology of subspace of $\Omega(G)/\sim$, we can write
$$V=(A/\sim)\cap (U/\sim)=(A\cap U)/\sim.$$
Here $U$ is an open subspace of $\Omega(G)$ with the property that
  for any $x\in U$, we have
$$\{y\in \Omega(G) \ : \ y\sim x\} \subset U.$$
The subspace $U\cap \Omega_{\a}(G)$ of $\Omega_{\a}(G)$ is  open in the direct limit topology
 (because the direct limit topology on
$\Omega_{\a}(G)$ is finer than the subspace topology, see Proposition
\ref{fine}) and non-empty
(because $\Omega_{\a}(G)$ is dense in $\Omega(G)$, see 
{\cite[section 3.5]{Pr-Se}}). For any 
$x\in U\cap \Omega_{\a}(G)$ we have
\begin{equation}\label{yyyy}(T^{\bC}\times \bC^*)x=\{y\in \Omega_{\a}(G) \ : \ y\sim x\}
 \subset U\cap \Omega_{\a}(G),
\end{equation}
which follows from Proposition \ref{propfirst} (a).
There exists $\lambda \in \check{T}$  such that
$\overline{\C}\cap U \neq \emptyset$ and  $\mu^{-1}(a)\cap \overline{\C}\neq \emptyset$.
Indeed,
we can pick $x\in \Omega_{\a}(G)\cap U$ (the intersection is non-empty, see above) and 
$y\in \Omega_{\a}(G)\cap \mu^{-1}(a)$ (the intersection is non-empty, since $a\in \mu(\Omega(G))=
\mu(\Omega_{\a}(G))$); by Proposition \ref{last}, there exists $\lambda \in \check{T}$ 
such that 
$x$ and $y$ are both  in $\overline{C_{\lambda}}$.

\noindent{\it Claim.} If $\lambda \in \check{T}$ is as above, then  $\overline{\C}^{ss}$
is a dense subspace of $\overline{\C}$ (here
$\overline{\C}$ is equipped with the direct limit topology it inherits from
$\Omega_{\a}(G)$).

To prove the claim, we consider again a Grassmannian ${\mathcal G}_n$
which contains  
$\overline{\C}$ as a $T^\bC \times \bC^*$-invariant closed subvariety. By the main theorem of
\cite{He-Mi}, there exists on $\mathcal{G}_n$ a  $T^{\bC}\times \bC^*$-invariant 
 very ample line bundle $L$ such that $ \mathcal{G}_n^{ss}=\mathcal{G}_n^{ss}(L)$.
 The latter space consists of all $L$-semistable points in 
 ${\mathcal G}_n$, that is points $x\in {\mathcal G}_n$
 such that there exists $k\ge 1$ and $s : X \to L^{\otimes k}$  
 equivariant    holomorphic  section       {\rm  with}  $s(x)\neq 0$
 (cf. e.g. \cite{Mu-Fo-Ki}).
  Consequently, $\mathcal{G}_n^{ss}$ is a Zariski open subspace
 of $\mathcal{G}_n$.
Since
$\overline{\C}^{ss}={\mathcal G}_n^{ss}\cap \overline{\C},$
 we deduce that $\overline{\C}^{ss}$ is a  Zariski open  subspace of 
 $\overline{\C}$. The space $\overline{\C}^{ss}$ is non-empty,
 since $\mu^{-1}(a)\cap  \overline{\C}\subset \overline{\C}^{ss}$.
Thus $\overline{\C}^{ss}$  is dense in $\overline{\C}$ with respect to the usual differential topology on the latter space: this can be deduced by using \cite[Theorem 2.33]{Mu}
for   $\overline{\C}$, which is an irreducible projective variety 
 (cf. \cite[p. 360]{Mi}).

From the claim we deduce that the intersection $ \overline{\C}^{ss}\cap U$
is non-empty (since $\overline{\C}\cap U$ is a non-empty subspace of 
$\overline{\C}$ which is open with respect to the direct limit topology).
By Lemma \ref{lema} we have
 $$\overline{\C}^{ss}\cap U=A\cap \overline{\C} \cap U,$$
  thus $$U \cap A \cap \Omega_{\a}(G)\neq \emptyset.$$
By equation (\ref{yyyy}), the quotient $(U \cap A \cap \Omega_{\a}(G))/ \sim$
is a (non-empty) subspace of $\Omega(G)/\sim$. It is contained in both
$V=(U\cap A)/\sim$ and $(A\cap \Omega_{\a}(G))/\sim$.  
Consequently, the intersection
$ V\cap [(A\cap 
\Omega_{\a}(G))/\sim]$ is non-empty.
  This finishes the proof.
  \hfill $\square$

\bibliographystyle{abbr}

\end{document}